\documentclass[11pt,a4paper]{article}
\usepackage[utf8]{inputenc}
\usepackage[english]{babel}
\usepackage{lmodern} 
\usepackage{amsmath}
\usepackage{amsfonts}
\usepackage{amssymb}
\usepackage{xcolor}

\usepackage{enumitem}

\newtheorem{Pp}{Property}
\newtheorem{C}[Pp]{Corollary}
\newtheorem{Po}[Pp]{Proposition}
\newtheorem{T}[Pp]{Theorem}
\newtheorem{Le}[Pp]{Lemma}

\def\QEDmark{\ensuremath{\square}\\}
\def\proof{\noindent\textbf{\textit{Proof. }}}
\def\endproof{\hfill\QEDmark}  

\usepackage{mathrsfs}	
\newcommand{\Stz}{\mathscr{S}}
\usepackage{aurical} 
\newcommand{\CB}{\text{\Fontauri B}_1}
\newcommand{\TFb}[1]{\widehat{b_{#1,n}}}

\author{Ir\`{e}ne Casseli
	\footnote{
			\href{mailto:irene.casseli@univ-amu.fr}{irene.casseli@univ-amu.fr}
		}
	}
\date{\small	\textit{Aix-Marseille Universit\'{e}, I2M UMR CNRS 7373, 39 Rue F. Joliot-Curie, 13453
		Marseille Cedex 13, France}} 
\title{Fixed points of the Berezin transform of polyanlytic Fock spaces}

\usepackage{url}
\usepackage{hyperref}
\hypersetup{
	breaklinks=true,
	colorlinks=true,                         
	linkcolor=red!60!black, 
	citecolor=cyan!90!black, 
	urlcolor=cyan!90!black,
	linkcolor=red!60!black,
	} 

\begin{document}

\maketitle

	\begin{abstract}
		We study the fixed points of the Berezin transform in polyanalytic Fock spaces of $\mathbb{C}$. We show that an $L^p$ function, $p\in[1,+\infty]$, with respect to the Lebesgue measure is invariant under this transformation if and only if it is harmonic. From this we deduce that the only bounded fixed points of the Berezin transform of polyanalytic Fock spaces are constant functions.
	\end{abstract}~\\
	
		
		


	\section{Introduction}
	
	The Berezin transform is a powerful tool in operator theory. Initiated by Berezin \cite{berezin_covariant_1972} in connection with quantum mechanic, the study of such a transformation provides a lot of information mainly in the field of Toeplitz, Hankel and composition operators. The interested reader will find many applications in the context of the classical Bergman space of the unit disk $A^2(\mathbb{D})$ in \cite{zhu_operator_2007}, and other in \cite{zhu_analysis_2012} for the classical Fock space $F^2(\mathbb{C})$ setting. 
	
	Let us make precise the link between the Berezin transform of a bounded function and Toeplitz operators. 
	In both situations of $H=A^2(\mathbb{D})$ and $H=F^2(\mathbb{C})$, $\mu$ will denote the relevant measure, $\kappa_z(w)$ the normalized reproducing kernel for $H$, and $P$ the orthogonal projection operator from $L^2(d\mu)$ onto $H$. Let $\langle f,g \rangle$ be the usual inner product of  $f,g\in L^2(d\mu)$. The Berezin transform of a bounded measurable function $f$ is given by the formula  
	\[ Bf(z) = \int f |\kappa_z|^2 d\mu =\langle f\kappa_z,\kappa_z \rangle.\]
	Because $\langle f\kappa_z,\kappa_z \rangle = \langle P(f\kappa_z),\kappa_z \rangle$ when $f$ is bounded, its Berezin transform  turns out to be intimately related to the densely defined Toeplitz operator $T_f\colon g\in H\mapsto P(fg)$. 
	
	The characterization of functions invariant under the Berezin transform first emerged in \cite{axler_commuting_1991} from the problem of determining commuting Toeplitz operators of the Bergman space with bounded harmonic symbols.
	In \cite{englis_functions_1994}, Engli\v{s} gave a general solution and we will follow his method in the present paper. 
	An other solution in the context of Bergman spaces 
	is due to Ahern, Flores and Rudin \cite{ahern_invariant_1993}. In the one dimentional case, the result is that as long as the Berezin transform of a function $f$ is well defined, $f$ is invariant if and only if it is harmonic. 
	
	For the classical Fock space $F^2(\mathbb{C})$, a Lipschitz estimate for the Berezin transform of a bounded function, together with a semigroup property, shows that  a bounded fixed point of the Berezin
	transform must be constant (see \cite{zhu_analysis_2012} for details). Here, we are interested in providing a necessary and sufficient condition for an $L^p$ function with respect to the Lebesgue measure to be invariant under the Berezin transform of polyanalytic Fock spaces of $\mathbb{C}$. We also give a generalization of the former result about bounded fixed points in the new setting as a consequence of our main theorem.\\
	
	Let us introduce polyanalytic functions as defined in \cite{balk_polyanalytic_1991}. For $n\in \mathbb{N}^*$, a function $f$ is called $n$-analytic on $\mathbb{C}$ if it satisfies the condition \[ \frac{\partial^{n}f}{\partial\overline{z}^{n}}=0\] in the whole complex plane. 
	
	We consider the Gaussian probability measure \[ d\mu(z)= \frac{\mathrm{e}^{-|z|^2}}{\pi} d\lambda(z) \] where $\lambda$ is the Lebesgue area measure on $\mathbb{C}$. 
	For $f$ and $g$, two square integrable functions with respect to $\mu$, $\langle f,g\rangle$ will denote the usual $L^2(\mathbb{C},d\mu)$ inner product.
	For $n\in\mathbb{N}^*$, the $n$-analytic Fock space $F^2_{n}(\mathbb{C})$ is the closed subspace in $L^2(\mathbb{C},d\mu)$ consisting of all square integrable $n$-analytic functions. The reproducing kernel of the Hilbert space $F^2_{n}(\mathbb{C})$ is given  by 
	\begin{equation}\label{Def_Noyau} K_{n}(z,w) = 
	L^1_{n-1}(|z-w|^2)\mathrm{e}^{ z\overline{w}}\end{equation} for all $z,~w\in\mathbb{C}$ (see \cite{abreu_structure_2010} or \cite{askour_espaces_1997} for instance) where $L^\beta_k$ is the generalized Laguerre polynomials\begin{equation}\label{Def_Laguerre} L^\beta_k (x)= \sum_{j=0}^{k}  \begin{pmatrix}
	k+\beta\\k-j
	\end{pmatrix} \frac{(-x)^j}{j!}.\end{equation} For $z\in\mathbb{C}$, we also introduce the normalized kernel function 
	\[k_z^{n}\colon w \in\mathbb{C} \mapsto \frac{K_{n}(w,z)}{\sqrt{K_{n}(z,z)}}.\]
	
	The Berezin transform of $f\in L^p(\mathbb{C},d\lambda)$ is  defined on $\mathbb{C}$ by the formula     \[ B_{n}{f} (z)= \langle fk_z^{n},k_z^{n}\rangle,\quad z\in \mathbb{C}. \]

	Here is the main result. 
	
	\begin{T}\label{Th_Points_fixes}
		Let $f\in L^p(\mathbb{C},d\lambda)$. Then $B_{n}f=f$ if and only if $f$ is harmonic on $\mathbb{C}$.
	\end{T}
	
	As a consequence we have the following corollary.

	\begin{C}\label{Th_Points_fixes_bornes=cst}
		Let $f\in L^\infty(\mathbb{C})$. The following conditions are equivalent:
		\renewcommand{\theenumi}{\textit{(\roman{enumi})}}
		\renewcommand{\labelenumi}{\textit{(\roman{enumi})}}
		\begin{enumerate}
			\item $B_{n}f=f$ ; \label{1}
			\item $f$ is harmonic on $\mathbb{C}$; \label{2}
			\item $f$ is constant. \label{3}
		\end{enumerate}
	\end{C}
	
	The equivalence of \ref{2} and \ref{3} directly follows  from Liouville's theorem. For $n=1$ this corollary is reduced to Proposition 3.27 in \cite{zhu_analysis_2012}.

	\section{Some preparatory results}
	
	In this section we discuss some preliminaries on which our proof is based.
	
	We begin by recalling the link connecting Laguerre polynomials and Bessel functions. Let $J_0$ the Bessel function of the first kind (see \cite{petiau_theorie_1955}) expressed by the integral representation on $\mathbb{R}$ : \[ J_0(x)= \frac{1}{\pi}\int_{0}^{\pi} \cos(x\cos\theta)d\theta.\]
	Using this definition, it is quite immediate to check the next properties. 
	
	\begin{Pp}~\label{Pp_J0}
		The function $J_0$ is even and satisfies $J_0(0)=1$ and $|J_0|<1$ on $\mathbb{R}^*$.
	\end{Pp}
	
	The following well known proposition (see \cite{szego_orthogonal_1975} p.103, Theorem 5.4) exhibits the strong interaction between classical Laguerre polynomials and $J_0$.
	
	\begin{Po} \label{Po_Bessel_Lag}
		For all $k\in\mathbb{N}$
		\begin{equation} 
		k! L_{k}^{0}(x)\mathrm{e}^{-x}= \int_0^{+\infty} \mathrm{e}^{-t} t^{k} J_{0 }(2\sqrt{xt})dt.
		\end{equation}
	\end{Po}~\\
		
	Now, for a precise statement of our results, we introduce some classical notations. As usual $\mathbb{C}$ is identified with $\mathbb{R}^{2}$. The space of infinitely differentiable functions on a domain $\Omega$ in $\mathbb{C}$ is denoted by $\mathcal{C}^\infty(\Omega)$, and its subspace of compactly supported functions by $\mathcal{D}(\Omega)$. For $x=(x_1,x_{2})\in\mathbb{R}^{2}$, and $\alpha =(\alpha_1,\alpha_{2})\in\mathbb{N}^{2}$
	a multi-index, we set $x^\alpha=x_1^{\alpha_1}x_{2}^{\alpha_{2}}$. If $j=1,2$, we denote by $\partial_j$
	the partial derivative with respect to
	the $j$th real variable $x_j$, and by $\partial^\alpha$ the differential 
	written $\partial^\alpha =\partial_1^{\alpha_1}\partial_{2}^{\alpha_{2}}$. As customary, $\Stz$ stands 
	for  the Schwartz class of rapidly decreasing functions on $\mathbb{R}^{2}$ that is the set of  $f\in\mathcal{C}^\infty(\mathbb{R}^{2})$ satisfying \[ \sup\limits_{x\in\mathbb{R}^{2}} |x^\alpha \partial^\beta f(x)|<+\infty\]for all multi-indices $\alpha,\beta\in\mathbb{N}^2$.
	We shall
	remind the reader
	that $\Stz$ is closed under the convolution, and thus is an algebra. Define the Fourier transform of $f\in L^1(\mathbb{C},d\lambda)$ by setting
	\[ \widehat{f}(z) = \int_{\mathbb{C}} f(w) \mathrm{e}^{-\mathrm{i}\text{Re}(z\overline{w})}d\lambda(w). \]
	We also record the following important fact. The Fourier transform is a homeomorphism from $\Stz$ onto
	itself.
	 
	As in \cite{englis_functions_1994}, the fact that a fixed point is harmonic will be a consequence of the following lemma (see \cite{brelot_ements_1959}).
	
	\begin{Le}[Weyl's lemma] \label{Le_Weyl}
		Let $\Omega$ be a domain in $\mathbb{R}^{2}$ and $u$ a  locally integrable function on $\Omega$ verifying 
		\[\int_{\mathbb{C}} u(z)\Delta \varphi(z)d\lambda(z) = 0
		\] 
		for all $\varphi\in \mathcal{D}(\Omega)$.
		Then
		$u \in \mathcal{C}^{\infty }(\Omega )$ and $u$ is harmonic. 
	\end{Le}

	\section{The Berezin transform}
	
	Even though we are interested in the Berezin transform of functions in $L^p(\mathbb{C},d\lambda)$, we will say a few words about the class of functions for which this transformation is well defined. 	
	This leads to define the probability measure 
	\[  d\nu(z)=\frac{{L^1_{n-1}(|z|^2)}^2}{n} d\mu(z) = \frac{{L^1_{n-1}(|z|^2)}^2\mathrm{e}^{-|z|^2}}{n\pi}  d\lambda(z) \]
	on the complex plane $\mathbb{C}$. We also introduce translation and translated reflexion maps of $\mathbb{C}$ as
	follows \[t_a\colon z\in\mathbb{C} \mapsto z+a,~~\phi_a\colon z\in\mathbb{C}\mapsto a-z,\] for each $a\in\mathbb{C}$. Given $f$ a Lebesgue measurable function on $\mathbb{C}$, the following conditions are equivalent for all $a\in\mathbb{C}$ : 
	\begin{eqnarray} 
		\int_{\mathbb{C}} |f(z)| |K_{n}(z,a)|^2 d\mu(z) < +\infty  
		&\Leftrightarrow&  f \circ t_a \in L^1(\mathbb{C}, d\nu)\label{I1}\\ &\Leftrightarrow & f \circ \phi_a \in L^1(\mathbb{C}, d\nu).\nonumber
	\end{eqnarray} 
	
	We let $\CB$ denote the set of measurable functions $f$, defined on $\mathbb{C}$, satisfying (\ref{I1}) for all $a\in\mathbb{C}$. 
	The following inclusions are obvious :
	\[ L^1(\mathbb{C},d\lambda)\subset \CB\subset L^1(\mathbb{C},d\nu). \]
	Furthermore, a simple estimate of the kernel $K_{n}$ yields that for all $a\in \mathbb{C}$ and $1\leq p< +\infty$, we have $K_{n}(.,a)\in L^p(\mathbb{C},d\mu)$. Thus we claim that \[\Big(L^p(\mathbb{C},d\lambda)\subset\Big) L^p(\mathbb{C},d\mu)\subset\CB\] 
	for any $1< p\leq +\infty$.
	
	The previous equivalences (\ref{I1}) allow us to define the  Berezin transform of $f\in \CB$ as follows :
	\begin{equation} \label{Def_TB_B1}
	B_{n}f (z) = \langle f k_z^{n},k_z^{n}\rangle = \int_{\mathbb{C}} f \,|k_z^{n}|^2 d\mu = \int_{\mathbb{C}} f \circ t_z d\nu
	= \int_{\mathbb{C}} f \circ\phi_z d\nu,\quad z\in\mathbb{C} 
	\end{equation}
	where equalities follow from a change of variables. 
	
	To establish the proof of our main statement, we also express the action of $B_{n}$ on an element $f \in \CB$
	as a convolution over $\mathbb{R}^{2}\simeq \mathbb{C}$. Using the last equality above, we have
	\begin{equation} \label{E_convol}
	B_{n}f 
	=  f\ast b_{n}
	\end{equation}
	where \begin{equation}\label{Def_b_{n}}
	b_{n}(z)= \frac{\mathrm{e}^{-|z|^2}}{n\pi}{L^1_{n-1}(|z|^2)}^2\end{equation} for any $z\in\mathbb{C}$.
	
	Obviously, $b_{n}\in \Stz$. Recall that $\Stz$ is an algebra under convolution. Thus observe that the formula (\ref{E_convol}) implies that $\Stz$ invariant subspace of $B_{n}$.	\\
	
	By Minkowski's inequality, it can be easily checked that $B_{n}$ is a bounded operator on $L^p(\mathbb{C},d\lambda)$ for any $1\leq p\leq +\infty$. The arguments leading to Theorem \ref{Th_Points_fixes_bornes=cst} appeal to the following description of the adjoint $B_{n}$ on $L^p(\mathbb{C},d\lambda)$ under the duality $L^p(\mathbb{C},d\lambda)^*=L^q(\mathbb{C},d\lambda)$ for $1\leq p<+\infty$ and $1/p+1/q=1$. As the next lemma shows this adjoint is nothing other than the operator $B_{n}$ on $L^q(\mathbb{C},d\lambda)$.
	
	\begin{Le} \label{Le_dualite} Let $1 \leq  p, q \leq +\infty$ satisfy $1/p+1/q=1$. The for all $f\in L^p(\mathbb{C},d\lambda)$ and $g\in L^q(\mathbb{C},d\lambda)$ we have   \[\displaystyle\int_{\mathbb{C}} B_{n}f\,g \,d\lambda = \int_{\mathbb{C}} f\,B_{n}g \,d\lambda.\]		
	\end{Le}
	
	{\proof Let $f\in L^p(\mathbb{C},d\lambda)$ and $g\in L^q(\mathbb{C},d\lambda)$. Since the Lebesgue measure $\lambda$ is invariant under translations, we check that  \[  \int_{\mathbb{C}} \int_{\mathbb{C}} |f(z)b_{n}(w-z)g(w)|  d\lambda(z)d\lambda(w) \leq ||f||_{L^p(\mathbb{C},d\lambda)} ||b_{n}||_{L^1(\mathbb{C},d\lambda)} ||g||_{L^q(\mathbb{C},d\lambda)}\]
	and the right-hand term	is finite. Now, we mention that $b_{n}$ is radial. So from this and from Fubini's theorem, it follows  that 
		\begin{eqnarray}
		\int_{\mathbb{C}} B_{n}f(w)g(w) \,d\lambda(w) 
		&=& \int_{\mathbb{C}} \int_{\mathbb{C}} f(z)b_{n}(w-z) d\lambda(z) g(w)
		d\lambda(w) \nonumber\\
		&=& \int_{\mathbb{C}} f(z) \int_{\mathbb{C}} g(w)b_{n}(z-w)
		d\lambda(w)d\lambda(z) \nonumber\\
		&=& \int_{\mathbb{C}} f(z) B_{n}g(z)d\lambda(z). \nonumber
		\end{eqnarray} 
		This completes the proof.
		\endproof}\\

	Since the Berezin transform commutes with both $t_a$ and $\phi_a$, for any $a\in\mathbb{C}$, the next property, well known in the case of classical Bergman and Fock spaces, is true.
	
	\begin{Pp} \label{PpHarmoPtsFixes}
		Each harmonic function in $\CB$ is invariant under $B_{n}$. 
	\end{Pp}
	
	\proof Let $u$ be harmonic in $\mathbb{C}$. Note that due to the fact that $k^{n}_z$ is a unit vector in $F^2_{n}(\mathbb{C})$, we have 
	\begin{eqnarray}
	1 = \langle k^{n}_z,k^{n}_z\rangle = B_n\textbf{1}(z)= \frac{2}{n}\int_{0}^{+\infty} {L^1_{n-1}(r^2)}^2\mathrm{e}^{- r^2} rdr \nonumber
	\end{eqnarray}
	where $\textbf{1}$ denotes the constant function $1$,
	and so 
	\begin{eqnarray}
	\frac{2}{n} \int_{0}^{+\infty} {L^1_{n-1}(r^2)}^2 \mathrm{e}^{-r^2} rdr =1.\label{Eq_Int_L^1_{n-1}}
	\end{eqnarray}
	Combining this with the mean value property of $u$, we get 
	\[B_{n}u(0) = \int_{\mathbb{C}} u(w) d\nu(w) = \frac{2}{n}\int_{0}^{+\infty}  u(0) {L^1_{n-1}(r^2)}^2\mathrm{e}^{-r^2} rdr =u(0) .\] This provides the desired result since $B_{n}$ commutes with translations.
	\endproof
	
	As a consequence, the Berezin transform of a harmonic function (whenever it makes sense) is harmonic. 
	Moreover, Theorem \ref{Th_Points_fixes} says that the converse of Property \ref{PpHarmoPtsFixes} is also true for functions of $L^p(\mathbb{C},d\lambda)$. 
	
	Since polyharmonic functions generalize harmonic functions in the same way as polyanalytic functions generalize holomorphic functions (see \cite{aronszajn_polyharmonic_1983}), one may ask whether Property \ref{PpHarmoPtsFixes} still holds if the harmonic hypothesis is replaced by its generalization.
	\begin{Pp} \label{Pp_mharmo}
		There exist polyharmonic functions of order $m$ with $m\geq 2$ in $\CB$ which are not fixed by $B_{n}$.
	\end{Pp}
		
	\proof	Consider the non-harmonic function $f\colon z\in\mathbb{C}\mapsto |z|^2$, which is polyharmonic of order $2$.
	Applying Fubini's theorem, we have \[B_{n}f(0) = \frac{2}{n}  
	\int_0^{+\infty} r^2 {L^1_{n-1}(r^2)}^2 \mathrm{e}^{-r^2} r dr. \] 
	But the terms on the right had side do not vanish because of the continuity and the positivity of the integrand.  
	This implies that the property is true.
	\endproof
		
	\section{The Fourier transform of $b_{n}$}
	
	In order to prove our main statement we will need to study the Fourier transform of the function $b_{n}$ defined in (\ref{Def_b_{n}}). 
	For all $z\in\mathbb{C}$,  \[ \widehat{b_{n}}(z) = \int_{\mathbb{C}} b_{n}(w) \mathrm{e}^{-\mathrm{i}\text{Re}(z\overline{w})}d\lambda(w). \]

	\begin{T} \label{T_TF_b} For all $z\in\mathbb{C}$, we have
		$\widehat{b_{n}}(z)=Q_{n}\Big(\frac{|z|^2}{4}\Big)\mathrm{e}^{-\frac{|z|^2}{4}}$ where $Q_{n}$ is a real polynomial.
	\end{T}
	
	\proof Fix $z\in\mathbb{C}$ and use formula (\ref{Def_Laguerre}) for Laguerre polynomials to obtain
	\begin{eqnarray}
	\widehat{b_{n}}(z)  &=& \frac{1}{n\pi} \sum_{k,l=0}^{n-1} \begin{pmatrix}
	n\\k+1
	\end{pmatrix}\begin{pmatrix}
	n\\l+1
	\end{pmatrix} \frac{(-1)^{k+l}}{k!l!} 
	\int_{\mathbb{C}}|w|^{2(k+l)}\mathrm{e}^{-|w|^2}  \mathrm{e}^{-\mathrm{i}\text{Re}({z}\overline{w})} d\lambda(w)  \nonumber\\
	&=&  \frac{1}{n\pi} \sum_{k,l=0}^{n-1} \begin{pmatrix}
	n\\k+1
	\end{pmatrix}\begin{pmatrix}
	n\\l+1
	\end{pmatrix} \frac{1}{k!l!} \Delta^{k+l}\Bigg(
	\int_{\mathbb{C}}\mathrm{e}^{-|w|^2}  \mathrm{e}^{-\mathrm{i} \text{Re}({z}\overline{w})} d\lambda(w) \Bigg)  \nonumber
	\end{eqnarray}
	where $\Delta = 4\,\partial\overline{\partial}$ is the usual Laplacian on $\mathbb{C}$.
	Recognizing the Fourier transform of a Gaussian function, the last relation can be expressed as
	\begin{eqnarray}
	\widehat{b_{n}}(z)  
	&=&  \frac{1}{n} \sum_{k,l=0}^{n-1} \begin{pmatrix}
	n\\k+1
	\end{pmatrix}\begin{pmatrix}
	n\\l+1
	\end{pmatrix} \frac{1}{k!l!} \Delta^{k+l}\bigg( \mathrm{e}^{-\frac{|z|^2}{4}}  \bigg) \nonumber.
	\end{eqnarray}
	Some calculations show that \[ \Delta^{k}\bigg( \mathrm{e}^{-\frac{|z|^2}{4}}  \bigg) = (-1)^kk!L^0_k\bigg( \frac{|z|^2}{4}  \bigg) \mathrm{e}^{-\frac{|z|^2}{4}}  \]
	for all nonnegative integer $k$ with $L^0_k$ the Laguerre polynomial defined in (\ref{Def_Laguerre}).
	This formula can also be seen as a particular case of a more general result for Hermite type polynomials (see  \cite{intissar_spectral_2006} for example).
	Therefore the Fourier transform of $b_{n}$ at $z$ can be expressed as
	\begin{eqnarray}
	\widehat{b_{n}}(z)   =  u_n	\bigg(
	\frac{|z|^2}{4}  \bigg)
	\end{eqnarray}
	where $u_n(x) = Q_n(x)\mathrm{e}^{-x}$ and
	\begin{eqnarray} \label{Def_Qn}
	Q_n(x) =  \frac{1}{n} \sum_{k,l=0}^{n-1} \begin{pmatrix}
	n\\k+1
	\end{pmatrix}\begin{pmatrix}
	n\\l+1
	\end{pmatrix}
	\begin{pmatrix}
	k+l\\k
	\end{pmatrix}
	(-1)^{k+l}L_{k+l}^{0}(x).
	\end{eqnarray}This completes the proof of the theorem.
	\endproof

	Now, our focus is put on $Q_n$ and $u_n$ explicitly given in the previous proof. For this purpose we give an alternative formula of $u_n$. Inserting in (\ref{Def_Qn}) the expression of the Laguerre polynomial given by Proposition \ref{Po_Bessel_Lag}, we obtain
	\begin{eqnarray} \nonumber
	u_n(x) & = & \frac{1}{n} \sum_{k,l=0}^{n-1} \begin{pmatrix}
	n\\k+1
	\end{pmatrix}\begin{pmatrix}
	n\\l+1
	\end{pmatrix}
	\frac{(-1)^{k+l}}{k!l!} \int_0^{+\infty} \mathrm{e}^{-t} t^{k+l} J_{0}(2\sqrt{xt})dt\\
	&=& \frac{1}{n}{\tiny} \int_0^{+\infty} \mathrm{e}^{-t} \Big(L^{1}_{n-1}(t)\Big)^2 J_{0}(2\sqrt{xt})dt.\label{Eq_un}
	\end{eqnarray}
	
	\begin{Le} 
		$Q_n(0)=1$. \label{Le_Qn(0)}
	\end{Le}
	
	\proof By (\ref{Eq_un}) evaluated at $x=0$,  \[Q_n(0) = \frac{1}{n}{\tiny} \int_0^{+\infty} \mathrm{e}^{-t} \Big(L^{1}_{n-1}(t)\Big)^2 dt. \]		
	Formula (\ref{Eq_Int_L^1_{n-1}})  in the proof of Property \ref{PpHarmoPtsFixes}, together with the change of variables $t=r^2$ yield the result.		
	\endproof

	\begin{Le} For all $x>0$,
		$|u_n(x)|<1$. \label{Le_un<1}
	\end{Le}
	
	\proof Let $x>0$. As a consequence of Property \ref{Pp_J0} of the Bessel function $J_0$ we have that 
	\begin{eqnarray} 
	\mathrm{e}^{-t} \Big(L^{1}_{n-1}(t)\Big)^2 |J_{0}(2\sqrt{xt})| < \mathrm{e}^{-t} \Big(L^{1}_{n-1}(t)\Big)^2  \nonumber
	\end{eqnarray}
	except possibly for a finite number of $t$.
	Again, in view of (\ref{Eq_un}), we get
	\begin{eqnarray} 
	\bigg|\frac{1}{n}{\tiny} \int_0^{+\infty} \mathrm{e}^{-t} \Big(L^{1}_{n-1}(t)\Big)^2 J_{0}(2\sqrt{xt})dt  \bigg|
	&<& \frac{1}{n}{\tiny} \int_0^{+\infty} \mathrm{e}^{-t} \Big(L^{1}_{n-1}(t)\Big)^2 dt \nonumber 
	\end{eqnarray}
	and the
	proof is complete in view of (\ref{Eq_Int_L^1_{n-1}}).
	\endproof

	The end of this section is devoted to the available function $H_{n}$ defined on $\mathbb{C}\simeq \mathbb{R}^2$ by \[ H_{n}(z)= \frac{1-\widehat{b_n}(z)}{-|z|^2}.
	\]
	In view of the last two lemmas, the following property is easy to prove.
	
	\begin{Pp} \label{Pp_Hdn}
		$H_{n}$ is a non-vanishing function of $\mathcal{C}^\infty(\mathbb{R}^2)$.
	\end{Pp}
	
	We now prove that the multiplication by
	$1/H_{n}$ maps $\Stz$ into itself.
	
	\begin{Le} \label{Le_Schwrtz_stable}
		$1/H_{n}\Stz \subset \Stz$.
	\end{Le} 
	
	\proof
	Let $f\in\Stz$ and fix two multi-indexes $\alpha,\beta \in\mathbb{N}^{2}$. Multidimensional Leibniz's rule gives that	for each $x\in\mathbb{R}^2$,	
	\[ \bigg| x^\alpha \partial^{\beta}  \frac{f}{H_{n}} (x) \bigg|
	\leq \sum_{\gamma_1=0}^{\beta_1}\sum_{\gamma_{2}=0}^{\beta_{2}}\begin{pmatrix}
	\beta_1\\\gamma_1
	\end{pmatrix}  \begin{pmatrix}\beta_{2}\\\gamma_{2}
	\end{pmatrix}
	\big|x^\alpha \partial^{\beta-\gamma}f (x)\big|
	\big|\partial^{\gamma}(1/H_{n})(x)   \big|.
	\]
	But for any $\gamma\in\mathbb{N}^{2}$, an induction argument shows that $\partial^{\gamma}(1/H_{n})$ is a quotian of a polynomial of the distributional derivatives of $H_{n}$ divided by $H_{n}^{|\gamma|}$ where the notation $|\gamma|$ 
	means 
	the length $\gamma_1+\gamma_{2}$ of the multi-index. 
	Combining properties of $H_{n}$ (Property \ref{Pp_Hdn}) with the rapid decay at infinity of the Schwartz function $f$, we conclude that the function $x\mapsto x^\alpha \partial^{\beta}  ({f}/H_{n}) (x)$ is bounded on compacts sets.
	
	Argue similarly for the remaining unbounded region to obtain  
	\begin{eqnarray} \bigg| x^\alpha \partial^{\beta}  \frac{f}{H_{n}} (x) \bigg| &=&\bigg| x^\alpha \partial^{\beta}  \bigg(\frac{|.|^2f}{1-\widehat{b_{n}}}\bigg) (x) \bigg|\nonumber\\
	&\leq& \sum_{\gamma_1=0}^{\beta_1}\sum_{\gamma_{2}=0}^{\beta_{2}}\begin{pmatrix}
	\beta_1\\\gamma_1
	\end{pmatrix}  \begin{pmatrix}\beta_{2}\\\gamma_{2}
	\end{pmatrix}
	\big|x^\alpha \partial^{\beta-\gamma}(|.|^2f) (x)\big|
	\big|\partial^{\gamma}(1-\widehat{b_{n}})^{-1}(x)   \big|\nonumber
	\end{eqnarray} 
	where $\partial^{\gamma}(1-\widehat{b_{n}})^{-1}$ is a polynomial in $1-\widehat{b_{n}}$ and its distributional derivatives up to $\gamma$ times $(1-\widehat{b_{n}})^{-|\gamma|}$.
	
	Next, observe that $\widehat{b_{n}}$, and all its partial derivatives, are continuous functions of $z\in\mathbb{C}\simeq\mathbb{R}^{2}$ and have limit $0$ as $|z|\to\infty$. Consequently $\partial^{\gamma}(1-\widehat{b_{n}})$ and $(1-\widehat{b_{n}})^{-|\gamma|}$ are bounded at infinity for any $\gamma$. 
	Again, thanks to the rapid decay of $f$, it follows that $x\mapsto x^\alpha \partial^{\beta}  ({f}/H_{n}) (x)$ remains bounded at infinity. 
	The desired result follows at once.
	\endproof	
	
	The lemma above allows us to show as a corollary the last result of this part.
	
	\begin{C} \label{Le_inclusion}
		$\quad \Delta \Stz \subset  (Id-B_{n})\Stz$. 
	\end{C}
	
	\proof Suppose that $f\in\Stz$. It follows from a simple computation that \[\widehat{\Delta f}(z)=-|z|^2\widehat{f}(z) = (1-\widehat{b}_n(z))\frac{\widehat{f}(z)}{H_{n}(z)}. \]
	Since the Fourier transform is a bijective operator on $\Stz$, we can apply Lemma \ref{Le_Schwrtz_stable} to find a function $g\in\Stz$ such that $\widehat{\Delta f}= (1-\widehat{b}_n)\widehat{g} = \widehat{(Id-B_{n})(g)}$, using also that the Fourier transform translates between convolution and multiplication of functions. 
	Since $(Id-B_{n})\Stz\subset \Stz$, we obtain the desired result.
	\endproof

	\section{Proof of the main result}	
	
	We now proceed with the proof of Theorem \ref{Th_Points_fixes}. 
	
	Suppose that $f\in L^p(\mathbb{C},d\lambda)$ with $B_{n}f=f$. From lemma \ref{Le_dualite}, we have  
	\[ \int_{\mathbb{C}} f \, (B_{n}g) \,d\lambda = \int_{\mathbb{C}} (B_{n}f)\,g \,d\lambda  =\int_{\mathbb{C}} f\,g \,d\lambda   \] for all $g\in\Stz$.  Then \[ \int_{\mathbb{C}} f \, (Id-B_{n})(g) \,d\lambda  =0, \quad \forall\,g\in\Stz. \] Therefore  \[ \int_{\mathbb{C}} f \, \Delta g \,d\lambda  =0, \quad \forall\,g\in\Stz \]   follows directly from Lemma \ref{Le_inclusion}.
	Thus $f$ is a weak solution of Laplace's equation and the desired results follows by applying Weyl's lemma.
	
	The converse statement of the theorem is a consequence of Property \ref{PpHarmoPtsFixes}.

	\section{Comments}
	
	\begin{enumerate}[leftmargin=0cm]
		\item 
		One verifies without difficulty that we might have stated all our results for the weighted polyanalytic Fock spaces of polyanalytic functions on $\mathbb{C}$ which are square integrable with respect to the weighted Gaussian measure\[ d\mu(z)= \frac{\alpha}{\pi} \mathrm{e}^{-\alpha |z|^2} d\lambda(z) \]where $\alpha$ is a positive parameter. Indeed the operation that ascribes to each function $f$ its dilation $z\mapsto f(\sqrt{\alpha z})$ establishes an isomorphism between $F^2_n(\mathbb{C})$ and the weighted polyanalytic Fock spaces with parameter $\alpha$. Hence for the sake of convenience, we choose $\alpha=1$.

		\item The attempt to extend directly the proof of Proposition 3.27 in \cite{zhu_analysis_2012} to the polyanalytic case encounter some difficulties. 
		We have not been able to establish Lipschitz estimate as (3.13) in \cite{zhu_analysis_2012} because the polyanalytic kernel seems not to be favorable to a semigroup property.

		\item As mentioned in \cite{englis_functions_1994} the proof of Theorem \ref{Th_Points_fixes} extends from $L^p(\mathbb{C},d\lambda)$ to tempered distribution $\Stz'$. Although it is important to note that harmonic functions are fixed by the Berezin transform but more functions than the harmonic ones are fixed points. The example due to Engli\v{s} \cite{englis_functions_1994} also provide in our setting a non harmonic function which is invariant under the Berezin transform. Given $a\in\mathbb{C}$, pick $f_a\colon z\in\mathbb{C}
		\mapsto \mathrm{e}^{a\text{Re}z}$. This function satisfies $\Delta f_a = a^2 f_a$ ; in particular, it is not harmonic. However carrying out the calculation for all $z\in\mathbb{C}$, we have 
		\begin{eqnarray}
		B_{n}f_a(z)
		&=&  \frac{f_a(z)}{n}  \sum_{k,l=0}^{n-1} \begin{pmatrix}
		n\\k+1
		\end{pmatrix}\begin{pmatrix}
		n\\l+1
		\end{pmatrix} \frac{(-1)^{k+l}}{k!l!}  S_{k+l}(a) \, \nonumber 
		\end{eqnarray}
		where 
		\begin{eqnarray}
		S_m(a)
		=   \sum_{j=0}^{m} \begin{pmatrix}
		m\\j
		\end{pmatrix} I_{2(m-j)}(0) I_{2j}(a) \, \nonumber 
		\end{eqnarray}
		for a nonnegative integer $m$ and with $I_k(a) =  \sqrt{\pi}^{-1}\int_\mathbb{R} x^k \mathrm{e}^{ax- x^2}dx$ when $k\in\mathbb{N}$. Proceeding by induction it is not difficult to obtain that there exist polynomials $P_k$ with the same parity as $k$ such that $I_k(a) = P_k(a)\mathrm{e}^{a^2/4}$. Thus one can find a polynomial $P$ satisfying $B_{n}f_a=P(a^2)\mathrm{e}^{{a^2}/4}f_a$. Now owing to \cite{chabat_introduction1_1990}
		pp. 284-285, $\mathrm{e}^{-z/4}-P(z)$ has infinitely many zeros. Let us choose $a$ to be a square root of a nonzero solution of $\mathrm{e}^{-z/4}=P(z)$ ; this enables us to conclude that $f_a$ is invariant under $B_{n}$.
		
		\item The multidimensional case seems a little bit different and we are not able to use the approach of this article in general. Namely as a consequence of a recent work of Hachadi and Youssfi (see \cite{hachadi_polyanalytic_2019}), the reproducing kernel of the $n$-analytic Fock space $F^2_{n}(\mathbb{C}^d)$, $d\in\mathbb{N}^*$, is given by 
		\begin{equation}\label{Def_Noyau_d} K_{d,n}(z,w) = 
		\prod_{j=1}^{d}L^1_{n-1}(|z_j-w_j|^2)\mathrm{e}^{ z_j\overline{w_j}}\end{equation} for all $z=(z_1,\dots,z_d)$, $w=(w_1,\dots,w_d)\in\mathbb{C}^d$. 
		We also define the Berezin transform $B_{d,n}$ as a convolution over $\mathbb{R}^{2d}\simeq \mathbb{C}^{d}$ as soon as it is possible by 
		\begin{equation} 
		B_{d,n}f 
		=  f\ast b_{d,n}
		\end{equation} 
		where \begin{equation}
		b_{d,n}(z)= \frac{\mathrm{e}^{-|z|^2}}{(n\pi)^d} \prod_{j=1}^{d}{L^1_{n-1}(|z_j|^2)}^2\end{equation} for any $z\in\mathbb{C}^d$. One can check that the Fourier transform of $b_{d,n}$ can be expressed for all $z=(z_1,\dots,z_n)\in\mathbb{C}^d$ as
		$\widehat{b_{d,n}}(z)=\mathrm{e}^{-\frac{|z|^2}{4}}\prod_{j=1}^{d}Q_{n}\big(\frac{|z_j|^2}{4}\big)$ where $Q_n$ is the polynomial of Theorem \ref{T_TF_b}. But as soon as $n\geq 2$ and $d\geq 2$, the function $H_{d,n}$ defined on $\mathbb{C}^d$ by \[ H_{d,n}(z)= \frac{1-\TFb{d}(z)}{-|z|^2} 
		\] is not smooth at the origin so that Lemma  \ref{Le_Schwrtz_stable} fails. Actually, one can show that $ \Delta \Stz \subset  (Id-B_{d,n})L^q$ for each $q\geq 2$. Therefore, the following turns to be true. 
		\begin{T}
			Let $1\leq p \leq 2$. If $f\in L^p(\mathbb{C}^d,d\lambda_d)$ satisfies $B_{d,n}f=f$ then necessarily $f$ is harmonic on $\mathbb{C}^d$.
		\end{T} In view of the kernel formula (\ref{Def_Noyau}), we do not know if the converse statement is also true in general. To finish, we make one additional remark regarding Property \ref{PpHarmoPtsFixes}. The one dimensional case easily generalizes as follows.
		\begin{Pp}
			Complex-valued pluriharmonic functions are invariant under $B_{d,n}$ whenever it makes sence. 
		\end{Pp}
		
	\end{enumerate}
	
	\section*{Acknowledgements}

	I am grateful to my doctoral advisors Professors S. Rigat and E. H. Youssfi for helpful discussions. I also thank Nizar Demni for providing an additional argument and  appropriate  references.

	\bibliographystyle{final} 

\end{document}